\documentclass[11pt]{article}
\usepackage[margin=1in]{geometry}
\usepackage{geometry}
\usepackage{amsmath}
\usepackage{graphicx}
\usepackage{authblk}
\usepackage{float}
\usepackage{verbatim}  
\usepackage{mathtools} 
\usepackage{subcaption}

\usepackage{mathrsfs} 
\usepackage[dvipsnames]{xcolor}
\usepackage{amsfonts} 

\usepackage[toc,page]{appendix}
\usepackage{indentfirst}

\title{\textsc{Understand Slope Limiter - Graphically}}
\author{\small{Ling Zou} \\
        \small{Nuclear Science and Engineering Division, Argonne National Laboratory \\ 9700 S. Cass Ave, Lemont, IL 60439 (lzou@anl.gov)}}
\date{}

\begin{document}

\maketitle

\textsc{Abstract.} In this article, we illustrate how the concept of slope limiter can be interpreted graphically,
i.e., how the slope of reconstructed piecewise linear function is limited by four bounding lines that connect cell-averaged data and its neighboring cell-averaged data.
It is then conjectured that the same graphical rule can be generalized from uniform mesh to non-uniform mesh,
such that the high-resolution total variance diminishing (TVD) region of slope limiter for non-uniform meshes can be obtained.


\section{Introduction}
Slope limiter is central to high-resolution\footnote{These methods are commonly referred to as high-resolution instead of second-order methods, as they drop to first-order accuracy at local extrema.}
total variance diminishing (TVD) finite volume methods,
which are widely used in the numerical solving of hyperbolic partial different equations.
These methods can produce high-order spatial accuracy in regions with smooth solutions,
and can avoid spurious oscillations near discontinuities, such as shock waves.
In one dimension, this is realized by performing piecewise linear data reconstruction in each finite volume cell,
such that better-quality edge values can be obtained to compute numerical fluxes on the interface between neighboring cells.

Owing to the early contributions from van Leer \cite{vanLeer_1974}, Sweby \cite{Sweby_1984}, and many other researchers,
mathematical theories have been developed to form the concept of flux/slope limiter.
Among these earlier groundbreaking researches, Sweby's diagram has been widely used to illustrate the second-order TVD region of slope limiters.
These theories are mathematically rigorous, however they are not physically straightforward or intuitive,
especially to people with engineering background (like me) instead of mathematics background.
Inspired by the work of Berger et al. \cite{Berger_2005}, who provided a different perspective to examine the properties of slope limiters,
we found that the concept can be further interpreted in a more straightforward and physically intuitive way using graphics - the main idea discussed in this article.

In the following sections, we will at first have a very short introduction to the finite volume method to layout the context where the slope limiter concept will be discussed.
Following that, we will then discuss the concept of slope limiter, and how it can be interpreted graphically.
In the last section, we propose a conjecture that will extend high-resolution TVD slope limiters from uniform mesh to non-uniform mesh, based on which the second-order TVD region of slope limiter for non-uniform meshes can be obtained.

\section{Finite Volume Method}
Considering the following simple hyperbolic partial differential equation
\begin{equation}
\label{eqn:linear-adv-eqn}
  \begin{split}
    u_t + f_x &= 0, \quad\quad x \in (-\infty, \infty)\\
    f &= au,
  \end{split}
\end{equation}
where $u_t \equiv \partial u(x, t) / \partial t$, $f_x \equiv \partial f(x, t) / \partial x$,
and $a$ is a constant non-zero real number.
This equation is also called the linear advection equation, and can be rewritten in a simpler form as
\begin{equation}
  u_t + a u_x = 0,
\end{equation}
which could be interpolated as the advection of a tracer material, with concentration $u$, in a fluid field with flow speed $a$.

\begin{figure}[h]
  \centering
  \includegraphics[scale=0.8]{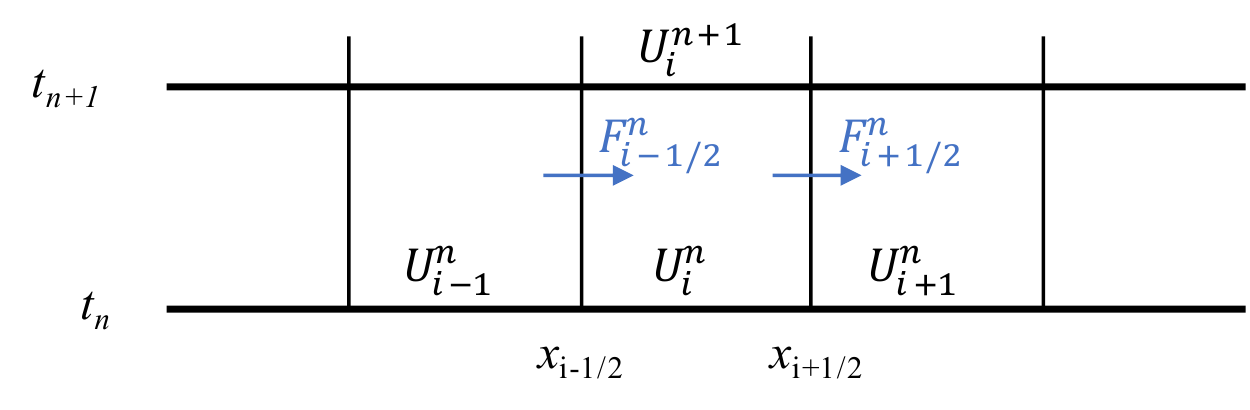}
  \caption{Illustration of a finite volume method in the $x$-$t$ space.}
  \label{fig:fvm}
\end{figure}

Let's solve the above hyperbolic equation \eqref{eqn:linear-adv-eqn} using the finite volume method on a mesh with a uniform cell size, as illustrated in figure \ref{fig:fvm}.
The i-th grid cell is denoted by
\begin{equation*}
  C_i = (x_{i-1/2}, x_{i+1/2}),
\end{equation*}
and we use $U_i^n$ to approximate the average value of $u$ over the i-th interval at time $t_n$
\begin{equation}
  U_i^n \approx \frac{1}{\Delta x} \int_{x_{i-1/2}}^{x_{i+1/2}} u(x, t_n) dx
  \equiv \frac{1}{\Delta x} \int_{C_i} u(x, t_n) dx
\end{equation}
in which, $\Delta x = x_{i+1/2} - x_{i-1/2}$, which is constant for the uniform-sized mesh.

Assuming that $U^n$ is already known, we wish to solve for $U^{n+1}$.
Consider the integral form of the conservation law, equation \eqref{eqn:linear-adv-eqn}, on the i-th interval,
\begin{equation}
  \frac{d}{dt} \int_{C_i} u(x, t) dx = f(u(x_{i-1/2}, t)) - f(u(x_{i+1/2}, t))
\end{equation}
which is obtained by applying Leibniz integral rule and the divergence theorem (Gauss's theorem).
Then, performing a time integral from $t_n$ to $t_{n+1}$ on the above equation, it is easy to see that
\begin{equation}
  \begin{split}
\frac{1}{\Delta x} \int_{C_i} u(x, t_{n+1}) dx = &\frac{1}{\Delta x} \int_{C_i} u(x, t_{n}) dx\\
&-\frac{1}{\Delta x} \left[\int_{t_n}^{t_{n+1}} f(u(x_{i+1/2}, t)) dt - \int_{t_n}^{t_{n+1}} f(u(x_{i-1/2}, t)) dt \right].
  \end{split}
\end{equation}
Now let's introduce a numerical flux, $F$, to approximate the average flux along the cell edge, i.e.,
\begin{equation}
  F_{i-1/2}^n \approx \frac{1}{\Delta t} \int_{t_n}^{t_{n+1}} f(u(x_{i-1/2}, t)) dt.
\end{equation}
If we can approximate this numerical flux based on the values $U^n$ (for explicit method),
we now have a fully discrete method to solve for $U^{n+1}$, i.e.,
\begin{equation}
  U_i^{n+1} = U_i^n - \frac{\Delta t}{\Delta x} \left(F_{i+1/2}^n - F_{i-1/2}^n\right).
\end{equation}

\subsection{REA Algorithm}

\begin{figure}
  \centering
  \begin{subfigure}[b]{0.3\textwidth}
    \centering
    \includegraphics[scale=0.8]{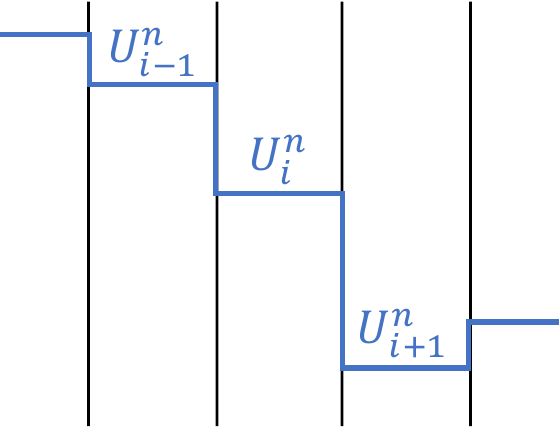}
    \caption{}
    \label{fig:reconstruct}
  \end{subfigure}
  \begin{subfigure}[b]{0.3\textwidth}
  \centering
  \includegraphics[scale=0.8]{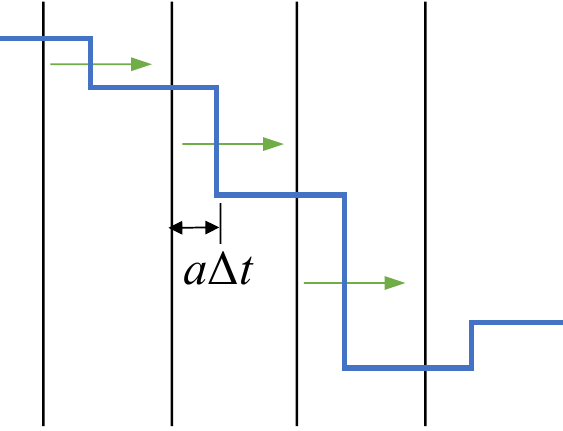}
  \caption{}
  \label{fig:evolve}
  \end{subfigure}
  \begin{subfigure}[b]{0.3\textwidth}
  \centering
  \includegraphics[scale=0.8]{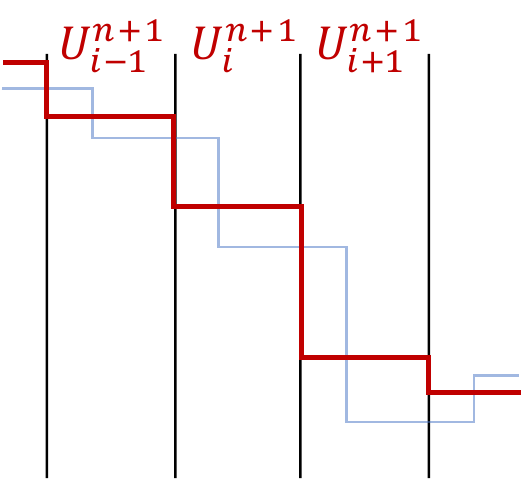}
  \caption{}
  \label{fig:average}
  \end{subfigure}
  \caption{(a) Reconstruct;
    (b) Evolve;
    (c) Average.}
  \label{fig:REA algorithm}
\end{figure}

For this advection equation,
the simplest way to compute the numerical flux is the so-called donor-cell upwind method.
Without loss of generality, assuming the flow is in the positive $x$ direction, i.e., $a > 0$.
For the donor-cell upwind method, the flux on the edges of each cell is entirely determined by its \textit{upwind} neighboring cell.
For example, on the left edge of the i-th cell, the flux is determined as
\begin{equation}
  F_{i-1/2}^n = a U_{i-1}^n,
\end{equation}
and similarly on the right edge,
\begin{equation}
  F_{i+1/2}^n = a U_i^n.
\end{equation}
Using this upwind method, the discrete equation becomes
\begin{equation}
  U_i^{n+1} = U_i^n - \frac{a \Delta t}{\Delta x} \left(U_i^n - U_{i-1}^n\right),
\end{equation}
from which $U_i^{n+1}$ can be updated.
This process can be illustrated as the process shown in figure \ref{fig:REA algorithm}.
We start from the known cell-averaged solutions, $U^n$, from time step $t_n$,
and then perform a simple piecewise constant reconstruction of the solution in each cell, shown as figure \ref{fig:reconstruct}.
We then evolve the solution using the upwind method within the time step, $\Delta t$,
i.e., the entire solution profile moves to the right (as $a > 0$) for a distance of $a \Delta t$, figure \ref{fig:evolve}.
Finally, perform local average in each cell to obtain the solution $U^{n+1}$ at $t_{n+1}$, shown as figure \ref{fig:average}.

This process can be seen as a simplified version of the reconstruct-evolve-average (REA) algorithm,
(see section 4.10 of LeVeque \cite{LeVeque_2002}).
We quote the REA algorithm and rewrite it as:

\noindent\textbf{REA Algorithm:}
\begin{enumerate}
  \item \textbf{Reconstruct} a piecewise polynomial function defined for all $x$ from the cell-averaged values.
        In the case demonstrated above, a picecwise constant function was reconstructed.
  \item \textbf{Evolve} the hyperbolic equation exactly (or approximately) with this reconstructed function for a time step $\Delta t$.
  \item \textbf{Average} the evolved function over each cell to obtain the new cell-averaged value for the new time step.
\end{enumerate}

\section{Slope Limiter}
In the previous section, in the context of finite volume method/REA algorithm,
we discussed the simplest reconstruction method, namely the piecewise constant reconstruction in each cell.
It is however well-understood that such an upwind scheme is only first-order.
To achieve higher-order (e.g. second-order) spatial accuracy, instead of piecewise constant reconstruction,
piecewise linear reconstruction is needed.

\begin{figure}[h]
  \centering
  \begin{subfigure}{0.6\textwidth}
  \centering
  \includegraphics[scale=0.8]{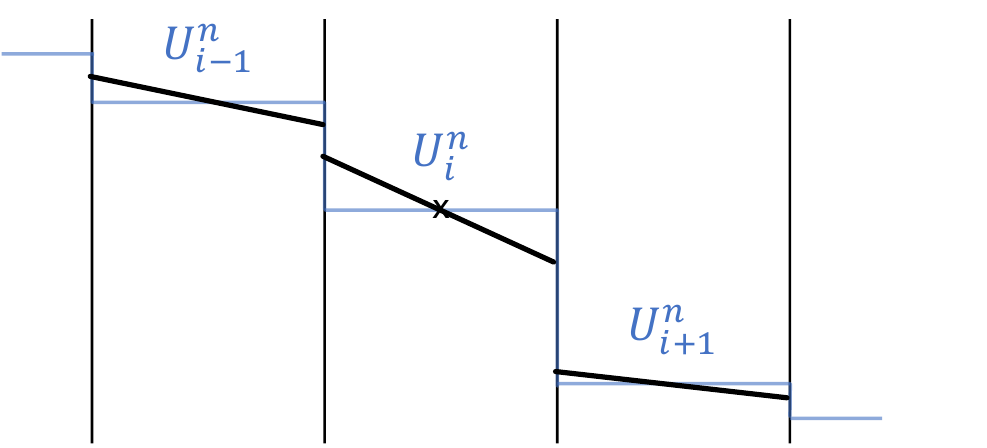}
  \caption{}
  \label{fig:linear-reconstruction}
  \end{subfigure}
  \begin{subfigure}{0.3\textwidth}
  \centering
  \includegraphics[scale=0.8]{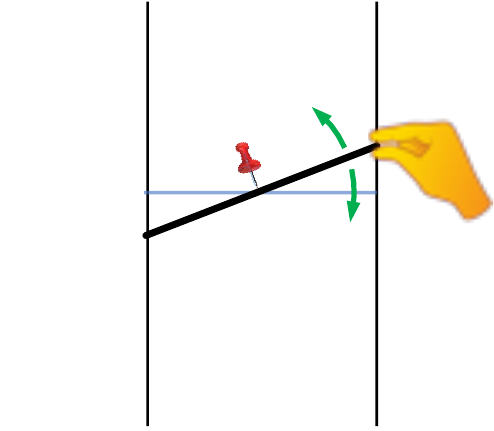}
  \caption{}
  \label{fig:linear-reconstruction-action}
  \end{subfigure}

  \caption{(a) Illustration of the piecewise linear reconstruction within each cell;
    (b) Illustration of how piecewise linear reconstruction is performed.}
\end{figure}

To maintain the total ``mass'' in each cell, e.g. the total amount of tracer material,
the reconstructed piecewise linear function must pass through, for example, the point ($x_i, U_i$) in the i-th cell,
marked as `x' in figure \ref{fig:linear-reconstruction}.
So this really leaves only one freedom we can play with, the slope of the piecewise linear function.
As illustrated in figure \ref{fig:linear-reconstruction-action},
to perform linear function reconstruction,
it is like that the center of the bar is pinned, while the bar could be rotated clockwise or counterclockwise (the `hand' symbol) about this pinned center.

If the local ``mass'' conservation is the only constraint, the slope can take any values.
However, we will see that the slope cannot take arbitrarily values, if we would like to construct a (high-resolution) TVD scheme.
For TVD schemes, the slope can only be chosen from a limited range, that is how \textit{slope limiter} got its name.
We do not repeat the concept of TVD here, which can be found in textbooks on finite volume method, such as LeVeque \cite{LeVeque_2002}.
In general, TVD method can be interpreted as a monotonicity-preserving method,
or more intuitively, a method that avoids non-physical overshoot/undershoot solutions.

\begin{figure}[H]
  \centering
  \begin{subfigure}{0.32\textwidth}
    \centering
    \includegraphics[scale=0.8]{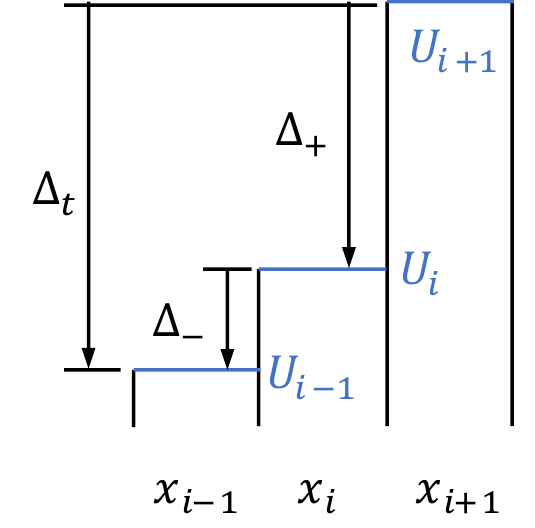}
    \caption{}
    \label{fig:deltas}
  \end{subfigure}
  \begin{subfigure}{0.32\textwidth}
    \centering
    \includegraphics[scale=0.8]{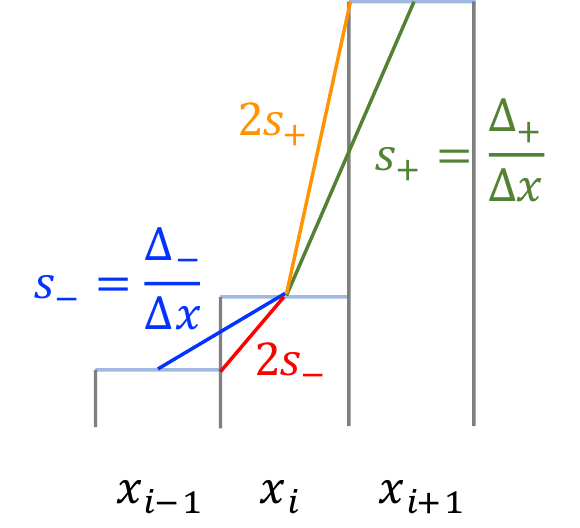}
    \caption{}
    \label{fig:s-minus-s-plus}
  \end{subfigure}
  \begin{subfigure}{0.32\textwidth}
    \centering
    \includegraphics[scale=0.8]{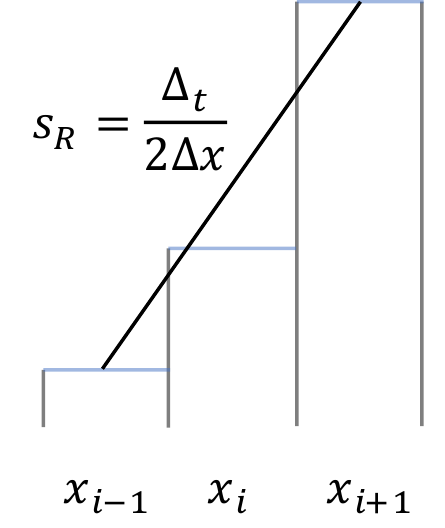}
    \caption{}
    \label{fig:reference-slope}
  \end{subfigure}

  \caption{(a) The differences between neighboring cell-averaged values;
    (b) Slopes relevant to construct slope limiter;
    (c) The reference slope;}
\end{figure}

Before we discuss the slope limiter, several important quantities ought to be introduced and discussed.
First, the differences between neighboring cell-averaged values are defined as:
\begin{equation}
  \begin{split}
    &\Delta_{-} \coloneqq U_i - U_{i-1} \\
    &\Delta_{+} \coloneqq U_{i+1} - U_i
  \end{split} \,,
\end{equation}
and
\begin{equation}
  \Delta_{t} = \Delta_{-} + \Delta_{+} \,,
\end{equation}
which are illustrated in figure \ref{fig:deltas}.
In above equations, we use $\coloneqq$ for `equal by definition'.
It is noted that their values can be positive or negative.
We then introduce the non-dimensional location indicator, $f$, defined as
\begin{equation}
  f \coloneqq \frac{\Delta_{-}}{\Delta_{t}}
\end{equation}
which indicates the relative location of $U_{i}$ in between $U_{i-1}$ and $U_{i+1}$.
It is noted that, if $U_{i-1}$, $U_{i}$, and $U_{i+1}$ are monotonic, either increasing or decreasing, we have $0 \leq f \leq 1$.
For all other conditions, $U_{i}$ is a local extremum, and we have $-\infty < f < 0$ or $1 < f < +\infty$.
It is also noted that, traditionally, the indicator is expressed as $r \coloneqq \Delta_{-} / \Delta_{+}$ (e.g., Sweby \cite{Sweby_1984}), which makes the concept of slope limiter less physically intuitive.

The two slopes associated with the `-' (left) and `+' (right) sides are defined as:
\begin{equation}
\begin{split}
  &{\color{blue}s_{-}} \coloneqq \frac{U_i - U_{i-1}}{\Delta x} = \frac{\Delta_{-}}{\Delta x}  \\
  &{\color{OliveGreen}s_{+}} \coloneqq \frac{U_{i+1} - U_{i}}{\Delta x} = \frac{\Delta_{+}}{\Delta x}
\end{split},
\end{equation}
which are illustrated in figure \ref{fig:s-minus-s-plus} as the {\color{blue}blue} and the {\color{OliveGreen}green} lines, respectively.
In the same figure, the {\color{red}red} line has a slope of {\color{red}$2s_{-}$}, and the {\color{orange}orange} line has a slope of {\color{orange}$2s_{+}$}.
We will see that, \textbf{these four bounding slopes (lines) play the utmost important roles to construct high-resolution TVD slope limiters}.

In figure \ref{fig:reference-slope}, a reference slope is also defined, which is the slope of the line that connects
the two neighboring data, i.e., ($x_{i-1}, U_{i-1}$) and ($x_{i+1}, U_{i+1}$),
\begin{equation}
  s_R \coloneqq \frac{U_{i+1} - U_{i-1}}{x_{i+1} - x_{i-1}} = \frac{\Delta_t}{2 \Delta x}.
\end{equation}

Using $s_R$, we can non-dimensionalize the slopes discussed above. It is easy to see:
\begin{equation}
  \begin{split}
    &{\color{blue}\phi_{-}} \coloneqq \frac{{\color{blue}s_{-}}}{s_R} = {\color{blue}2f} \,,     \\
    &{\color{OliveGreen}\phi_{+}} \coloneqq\frac{\color{OliveGreen}s_{+}}{s_R} = {\color{OliveGreen}2(1-f)}  \,,
  \end{split}
\end{equation}
in which $\color{blue}\phi_{-}$ and $\color{OliveGreen}\phi_{+}$ are the non-dimensional slopes of the {\color{blue}blue} and the {\color{OliveGreen}green} lines, respectively.
The {\color{red}red} line has a non-dimensional slope of $\color{red}4f$; and the {\color{orange}orange} line has a non-dimensional slope of $\color{orange}4(1-f)$.

\begin{figure}[h]
  \centering
  \begin{subfigure}{0.32\textwidth}
    \centering
    \includegraphics[scale=0.8]{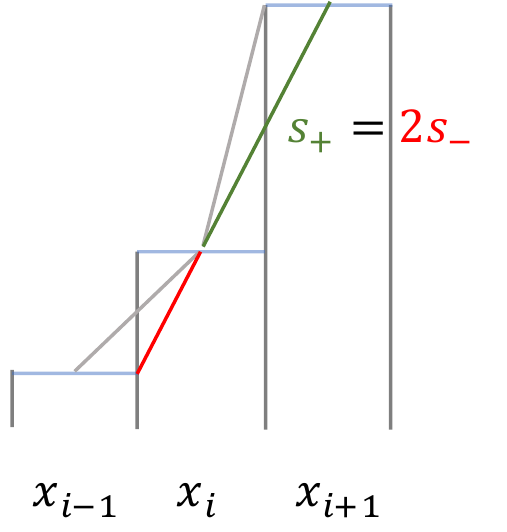}
    \caption{$s_{+} = 2 s_{-}$, when $f = 1/3$}
    \label{fig:f-1-third}
  \end{subfigure}
  \begin{subfigure}{0.32\textwidth}
    \centering
    \includegraphics[scale=0.8]{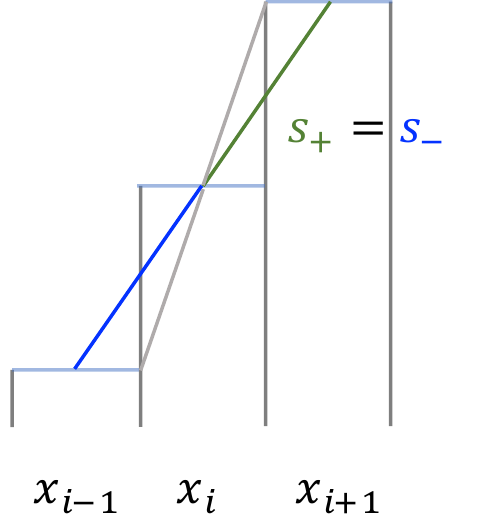}
    \caption{$s_{+} = s_{-}$, when $f = 1/2$}
    \label{fig:f-1-half}
  \end{subfigure}
  \begin{subfigure}{0.32\textwidth}
    \centering
    \includegraphics[scale=0.8]{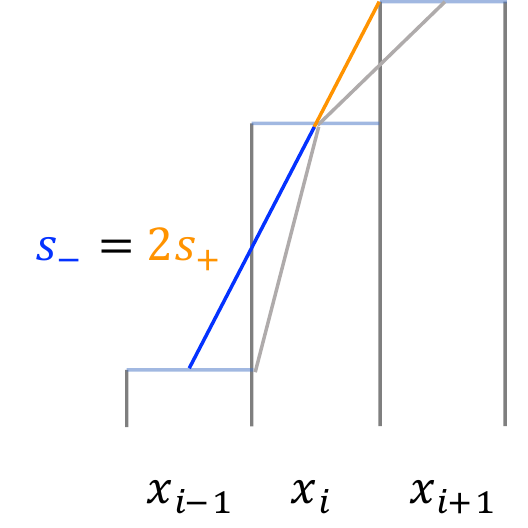}
    \caption{$s_{-} = 2 s_{+}$, when $f = 2/3$}
    \label{fig:f-2-third}
  \end{subfigure}
  \caption{Several special cases when two of the four slopes are equal. It is noted that cases (a) and (c) are symmetric about (b).}
  \label{fig:special-cases}
\end{figure}

There are several special cases when two of the four slopes are equal,
which are illustrated in figure \ref{fig:special-cases}.
The special values of $f$ will define the intervals where the slope limiter behaves differently.

Finally, let $s$ be the slope of the reconstructed piecewise linear function in the i-th cell.
We would like to learn how $s$ should behave when $U_i$ is at different relative locations in between $U_{i-1}$ and $U_{i+1}$.
Let's define the non-dimensional slope of the reconstructed piecewise linear function, $\phi$, as
\begin{equation}
  \phi \coloneqq \frac{s}{s_R}.
\end{equation}
It is now equivalent to ask: how $\phi$ should behave given different $f$.
We will learn that, to achieve high-resolution TVD scheme, $\phi$ can only be expressed as a nonlinear function of $f$, i.e., $\phi=\phi(f)$,
which will be dictated by the four colored lines shown in figure \ref{fig:s-minus-s-plus}.

\subsection{TVD Rule}

The rule for slope limiters to achieve TVD is very simple.
The mathematical theory behind the rules is well explained in Sweby \cite{Sweby_1984},
which we do not intend to repeat.
Instead, we will illustrate the two important rules in a graphical way.

\begin{figure}[h]
  \centering
  \begin{subfigure}[b]{0.3\textwidth}
    \centering
    \includegraphics[scale=0.8]{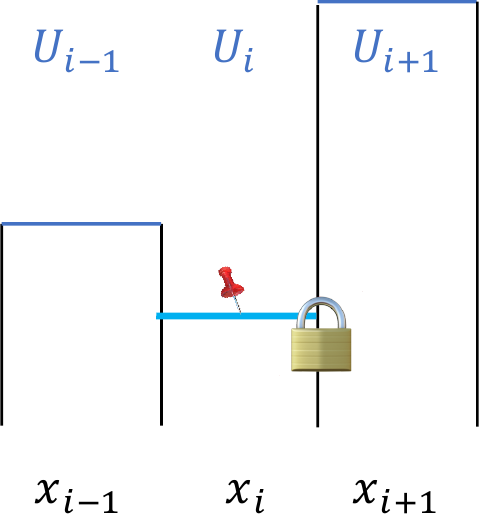}
    \caption{}
    \label{fig:local-extremum}
  \end{subfigure}
  \begin{subfigure}[b]{0.3\textwidth}
    \centering
    \includegraphics[scale=0.8]{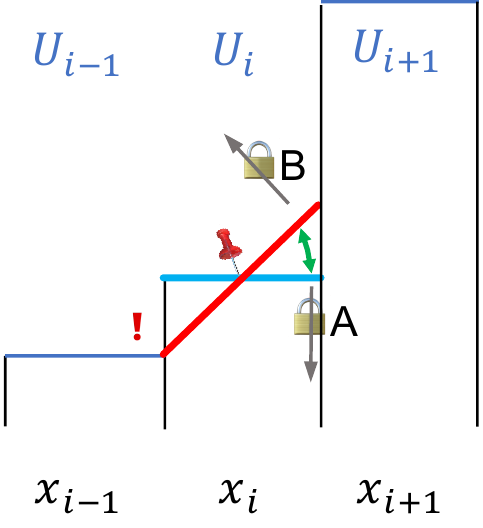}
    \caption{}
    \label{fig:tvd-1}
  \end{subfigure}
  \begin{subfigure}[b]{0.3\textwidth}
    \centering
    \includegraphics[scale=0.8]{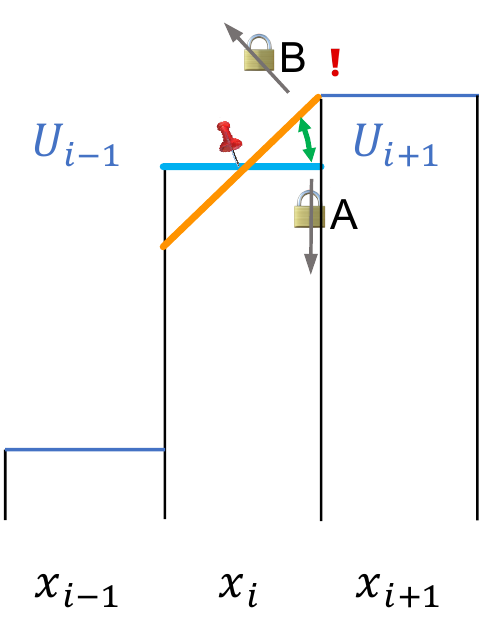}
    \caption{}
    \label{fig:tvd-2}
  \end{subfigure}

  \caption{
    The limiting
    (a) $U_i$ is a local extremum;
    (b) $U_i$ is in between $U_{i-1}$ and $U_{i+1}$, closer to $U_{i-1}$;
    (c) $U_i$ is in between $U_{i-1}$ and $U_{i+1}$, closer to $U_{i+1}$.}
\end{figure}

\noindent\textbf{TVD rules:}
\begin{enumerate}
  \item If $U_i$ is a local extremum, $-\infty < f < 0$ and $1 < f < +\infty$,
        the piecewise linear reconstruction must have a zero-slope, i.e., piecewise constant, as illustrated in figure \ref{fig:local-extremum}.

  \item For all other conditions ($0 \leq f \leq 1$), $U_i$ is in between $U_{i-1}$ and $U_{i+1}$.

        First, the reconstructed slope must not be negative (positive) for monotonically increasing (decreasing) data,
        as illustrated by the `lock'-A symbols in figures \ref{fig:tvd-1} and \ref{fig:tvd-2}.

        Second, the edge values of the linear function must not go beyond the neighboring cell-averaged values.
        This is illustrated in figure \ref{fig:tvd-1}, when $U_i$ is closer to $U_{i-1}$,
        the reconstructed linear function must not go beyond the {\color{red}red} line to allow its left-edge value to pass $U_{i-1}$.
        This is illustrated by the `!' symbol and the `lock'-B symbol in figure \ref{fig:tvd-1}.
        Similarly, figure \ref{fig:tvd-2} illustrates the case when $U_i$ is closer to $U_{i+1}$.

        In both figures \ref{fig:tvd-1} and \ref{fig:tvd-2}, the TVD region is simply bounded by the zero-slope blue line and the {\color{red}red}/{\color{orange}orange} line.
        Anything between them, marked by the green arrow, is acceptable as a TVD scheme.
\end{enumerate}

To summarize, for TVD schemes, the non-dimensional slope of the reconstructed piecewise linear function will respect the following conditions:
\begin{equation}
\left\{
  \begin{aligned}
    &0 \leq \phi(f) \leq 4f;      &\text{for} \quad 0 \leq f \leq 0.5      \\
    &0 \leq \phi(f) \leq 4(1-f);  &\text{for} \quad 0.5 < f \leq 1         \\
    &\phi(f) = 0;                 &\text{otherwise}
  \end{aligned}
\right. \,;
\end{equation}
or in a more compact form:
\begin{equation}
\left\{
  \begin{aligned}
    &0 \leq \phi(f) \leq \text{min}\left( 4f, 4(1-f)\right);      &\text{for} \quad 0 \leq f \leq 1      \\
    &\phi(f) = 0;                 &\text{otherwise}
  \end{aligned}
\right. \,.
\end{equation}
On the $\phi$-$f$ plot, the TVD region is illustrated as the shaded area in figure \ref{fig:tvd-phi-f}.

\subsection{High-resolution TVD}
In the previous section, we have discussed the slope limiter region for being TVD.
If we would like to further improve the slope limiter to achieve higher order spatial algorithm,
clearly, it will only be a subset of the TVD scheme.
The mathematical derivation of high-resolution TVD schemes can also be found in Sweby \cite{Sweby_1984}.
Again, we will proceed our discussion on a set of monotonically increasing data,
such that the four slopes, $s_{-}$, $2 s_{-}$, $s_{+}$, and $2 s_{+}$, are all positive.

The rule for high-resolution TVD is fascinatingly simple:

\vspace{0.1in}
\noindent\fbox{\parbox{\textwidth}{
    \textbf{High-resolution TVD slopes must lay in between the two smallest among the four bounding slopes.}
  }
}
\vspace{0.1in}

\begin{figure}[h]
  \centering
  \begin{subfigure}{0.18\textwidth}
    \centering
    \includegraphics[scale=0.7]{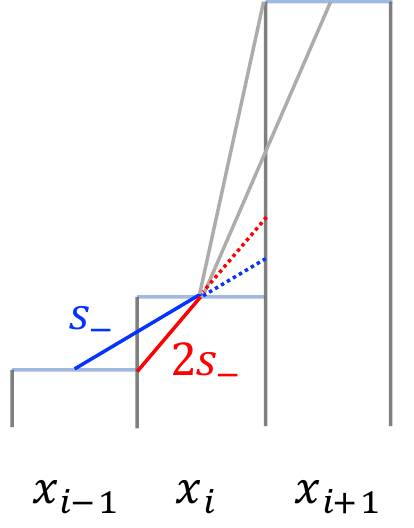}
    \caption{$0 < f < 1/3$}
    \label{fig:high-reso-int-1}
  \end{subfigure}
  \begin{subfigure}{0.18\textwidth}
    \centering
    \includegraphics[scale=0.7]{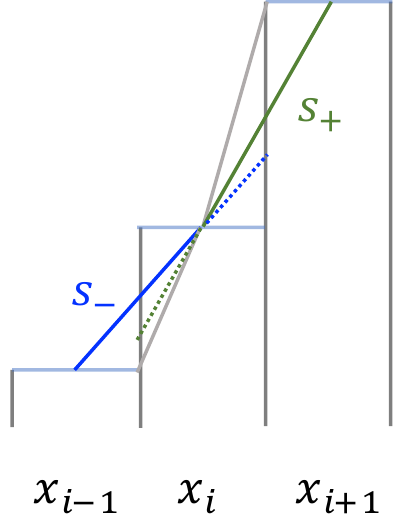}
    \caption{$1/3 < f < 1/2$}
    \label{fig:high-reso-int-2}
  \end{subfigure}
  \begin{subfigure}{0.18\textwidth}
    \centering
    \includegraphics[scale=0.7]{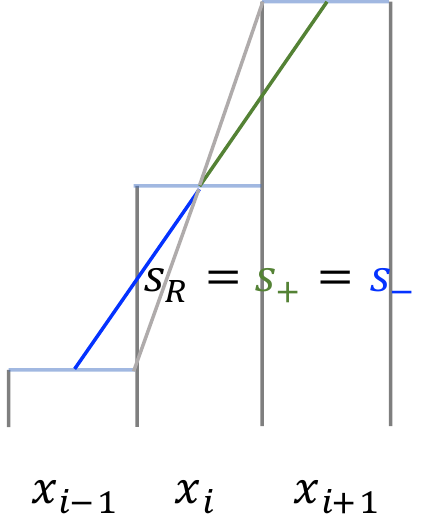}
    \caption{$f = 1/2$}
    \label{fig:high-reso-int-3}
  \end{subfigure}
    \begin{subfigure}{0.18\textwidth}
    \centering
    \includegraphics[scale=0.7]{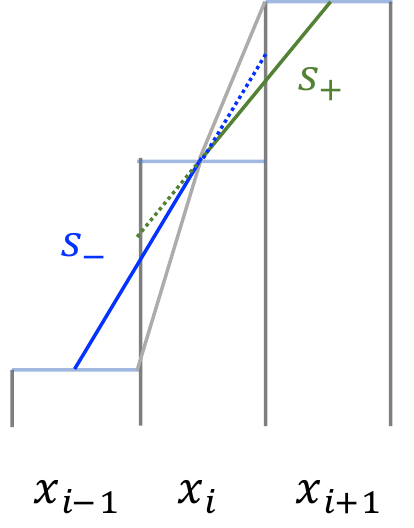}
    \caption{$1/2 < f < 2/3$}
    \label{fig:high-reso-int-4}
  \end{subfigure}
    \begin{subfigure}{0.18\textwidth}
    \centering
    \includegraphics[scale=0.7]{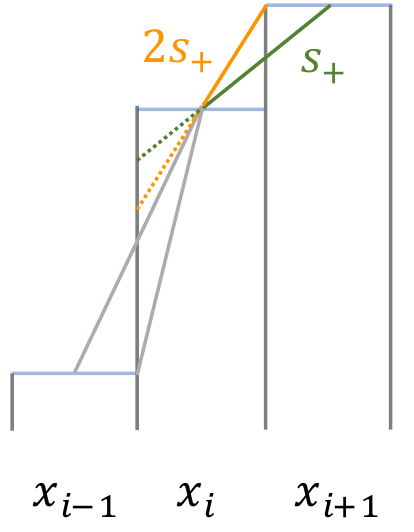}
    \caption{$2/3 < f < 1$}
    \label{fig:high-reso-int-5}
  \end{subfigure}

  \caption{Illustration of how slope should be limited to achieve high-order TVD for different conditions.
    For all cases, the acceptable reconstructed linear function must lay between the two colored lines, which have the smallest slopes among the four lines.
    For condition (c), the reconstructed linear function must have a slope of $s_R$.}
  \label{fig:tvd-intervals}
\end{figure}

\begin{figure}[h]
  \centering
  \begin{subfigure}{0.45\textwidth}
    \centering
    \includegraphics[scale=0.6]{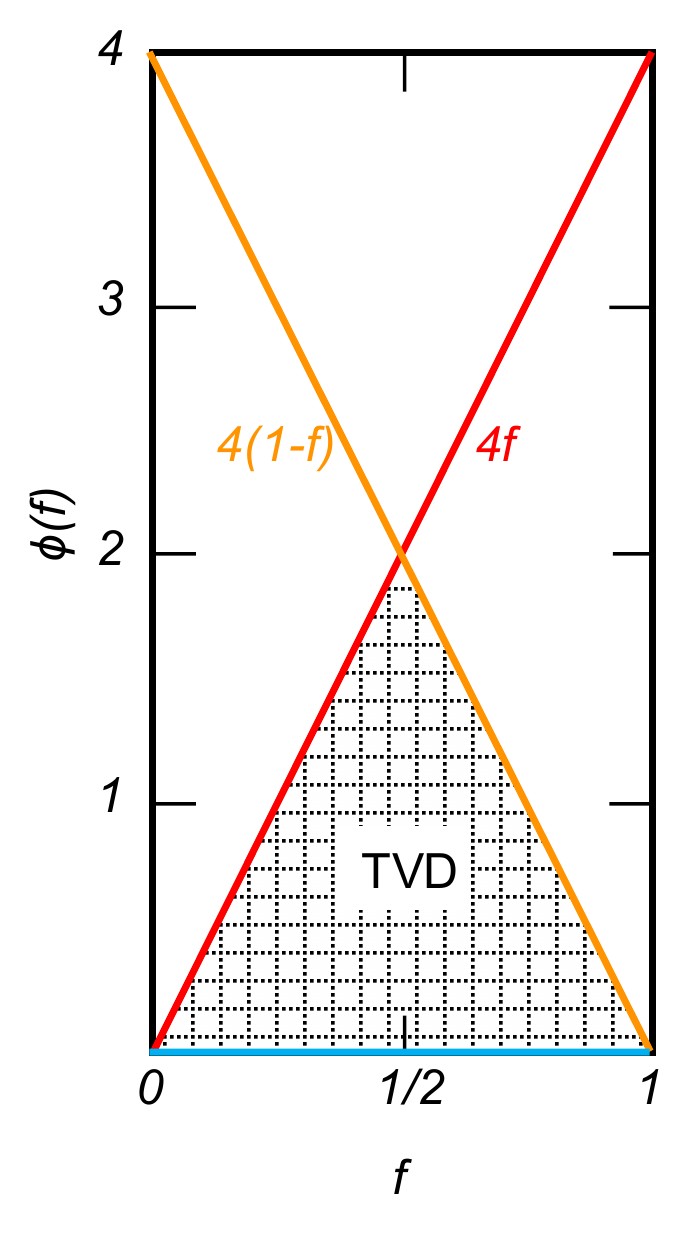}
    \caption{}
    \label{fig:tvd-phi-f}
  \end{subfigure}
  \begin{subfigure}{0.45\textwidth}
    \centering
    \includegraphics[scale=0.6]{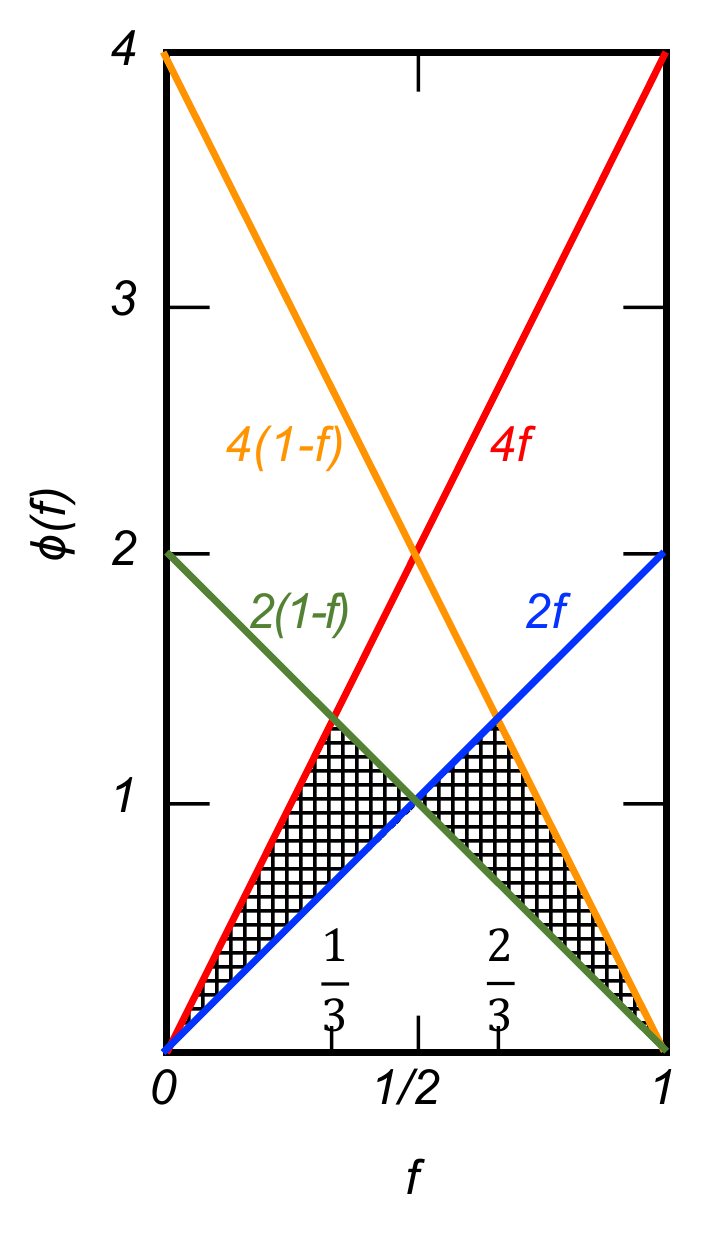}
    \caption{}
    \label{fig:high-resolution-tvd-phi-f}
  \end{subfigure}

  \caption{Illustration of (a) the TVD region (shaded); and (b) the high-resolution TVD region (shaded) on the $\phi$-$f$ plot.
    It is clear that the high-resolution TVD region is a subset of the TVD region.}
  \label{fig:phi-f-plot}
\end{figure}

Figure \ref{fig:tvd-intervals} illustrates the two smallest slopes for various $f$ values.
Especially, figure \ref{fig:high-reso-int-3} shows a special case where the two smallest slopes are identical, and also equal to the reference slope.
On the $\phi$-$f$ plot, the high-resolution TVD region can be easily found in between the two lines having smallest values among the four.
This is illustrated as the shaded area in figure \ref{fig:high-resolution-tvd-phi-f}.
It is noted that any high-resolution TVD schemes must pass the point:
\begin{equation}
  \phi\left(\frac{1}{2}\right) = 1,
\end{equation}
i.e., the special case shown in figure \ref{fig:high-reso-int-3}, which means that, if the solution is already a linear function, the reconstructed linear function must follow the same linear function to maintain the linearity.
This is the special case as illustrated in figure \ref{fig:high-reso-int-3}.

\begin{figure}[H]
  \centering
  \begin{subfigure}{0.3\textwidth}
    \centering
    \includegraphics[scale=0.8]{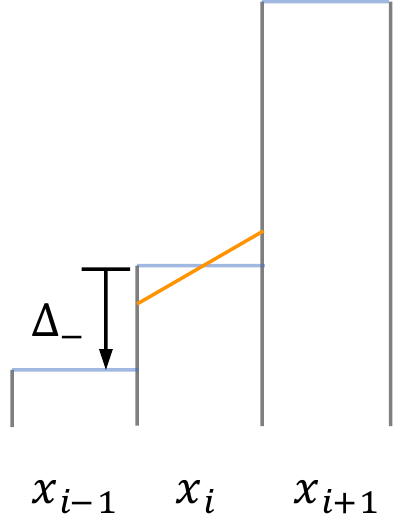}
    \caption{}
    \label{fig:sym-1}
  \end{subfigure}
  \begin{subfigure}{0.3\textwidth}
    \centering
    \includegraphics[scale=0.8]{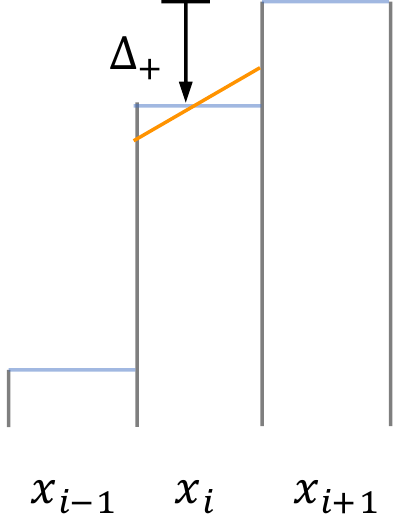}
    \caption{}
    \label{fig:sym-2}
  \end{subfigure}

  \caption{Illustration of a symmetric condition about $f=1/2$, $\Delta_{-}(a) = \Delta_{+}(b)$ and $\Delta_{t}(a) = \Delta_{t}(b)$.}
  \label{fig:symmetric}
\end{figure}

\begin{figure}[H]
  \centering
  \includegraphics[scale=0.4]{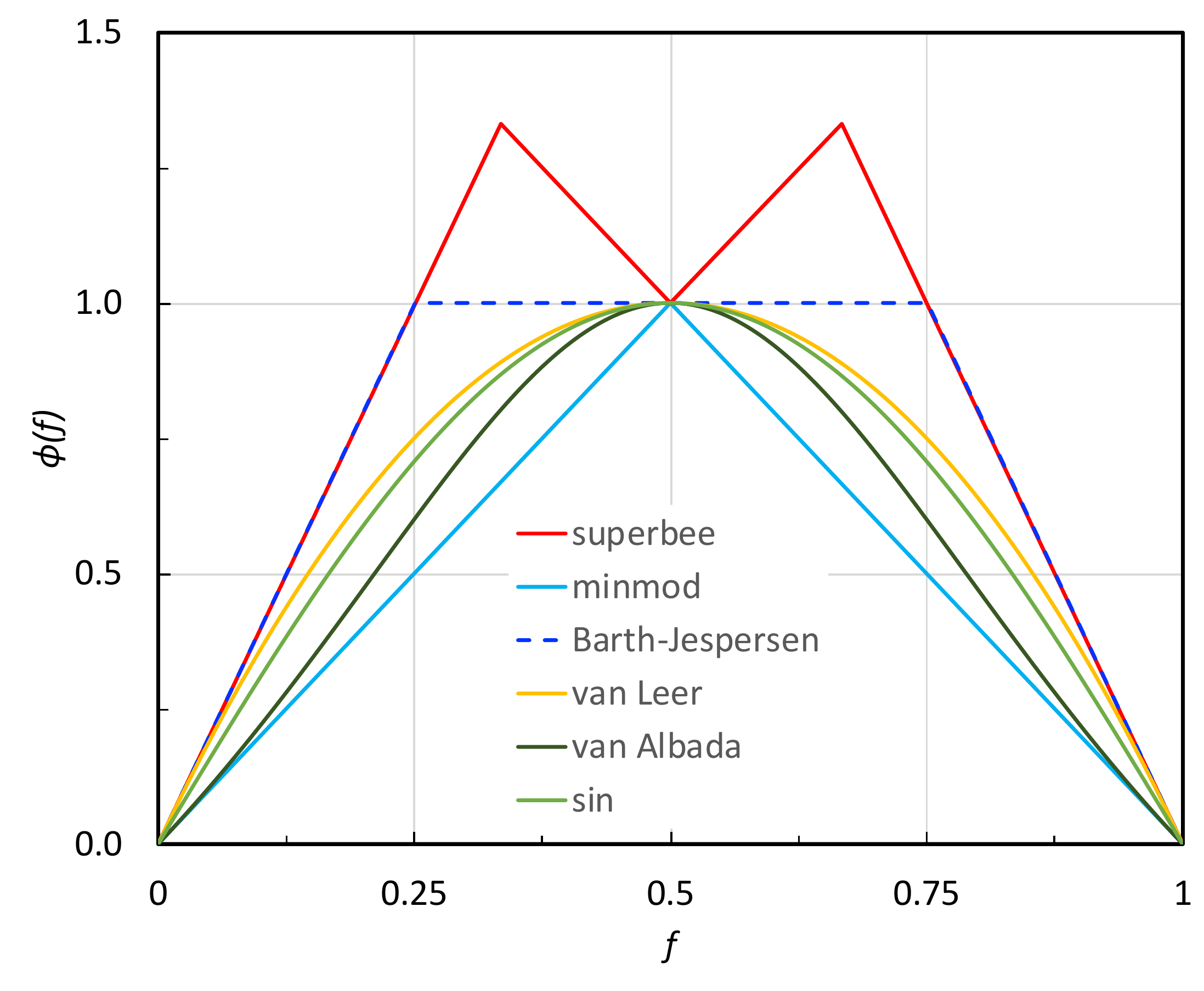}
  \caption{Illustration of several symmetric high-resolution TVD slope limiters.}
  \label{fig:slope-limiter-examples}
\end{figure}

\textbf{Symmetry:} It is noted that some high-resolution TVD scheme are so-called being symmetric.
This is can be illustrated in figure \ref{fig:symmetric}, in which conditions (a) and (b) are symmetric about $f = 1/2$,
and it seems natural to construct both linear functions with the same slope.
This requires the following condition:
\begin{equation}
  \phi(1-f) = \phi(f).
\end{equation}
Many well-known high-resolution TVD schemes do possess this symmetric feature.
Figure \ref{fig:slope-limiter-examples} shows several such examples,
which also include the \textit{sin} slope limiter proposed by Berger et al. \cite{Berger_2005} to demonstrate how a new symmetric slope limiter could be easily constructed.
However, we note that being symmetric is NOT required to be a high-resolution TVD scheme.

The mathematical form of these example limiters are given as follows.

\begin{itemize}
  \item \textbf{superbee}, the upper bound of the high-resolution TVD region
    \begin{equation}
    \phi(f) = \left\{
      \begin{aligned}
      &4f           & 0 \leq f \leq 1/3 \\
      &2(1-f)       & 1/3 < f \leq 1/2  \\
      &2f           & 1/2 < f \leq 2/3 \\
      &4(1-f)       & 2/3 < f \leq 1
      \end{aligned}
    \right.
    \end{equation}

  \item \textbf{minmod}, the lower bound of the high-resolution TVD region
  \begin{equation}
    \phi(f) = \text{min}[2f, 2(1-f)]
  \end{equation}

  \item \textbf{Barth-Jespersen}, also known as \textbf{MC}
  \begin{equation}
    \phi(f) = \text{min}[1, 4f, 4(1-f)]
  \end{equation}

  \item \textbf{van Leer}
  \begin{equation}
    \phi(f) = 4 f (1-f)
  \end{equation}

  \item \textbf{van Albada}
  \begin{equation}
    \phi(f) = \frac{2f(1-f)}{f^2 + (1-f)^2}
  \end{equation}

  \item \textbf{sin}
  \begin{equation}
    \phi(f) = \text{sin} (\pi f)
  \end{equation}
\end{itemize}

\section{Non-uniform Mesh}
To the author's best knowledge, there do not exist a complete theoretical derivation to define high-resolution TVD region on non-uniform meshes.
Berger et al. \cite{Berger_2005} provided a theoretical derivation to define the upper bound of TVD region and
the generalized \textit{min} slope limiter for high-resolution schemes on non-uniform meshes.
They also specified the necessary condition when dealing with linear data (see $f_2$ in equation \eqref{eqn:special_f_irregular_mesh}).

In the previous section, we have interpreted the high-resolution TVD region on uniform meshes graphically using the bounding slopes.
It seems natural that the same rule can be applied to irregular meshes,
which is proposed as the following conjecture:

\vspace{0.1in}
\noindent\fbox{\parbox{\textwidth}{
\textbf{Conjecture 1: \\
{}\\
Uniform mesh is a special case of the non-uniform mesh, and thus the high-resolution TVD rule for uniform meshes can be generalized and applied to non-uniform meshes.}
}
}
\vspace{0.1in}

To illustrate this, the reference slope and non-dimensional slopes need to be defined first.
Following Berger et al. \cite{Berger_2005}, we first define the mesh stretching ratios:
\begin{equation}
\begin{split}
  &a \coloneqq \frac{\Delta x_{i-1}}{\Delta x_{i}} \\
  &b \coloneqq \frac{\Delta x_{i+1}}{\Delta x_{i}}
\end{split}
\end{equation}

\begin{figure}[H]
  \centering
  \begin{subfigure}{0.45\textwidth}
    \centering
    \includegraphics[scale=0.8]{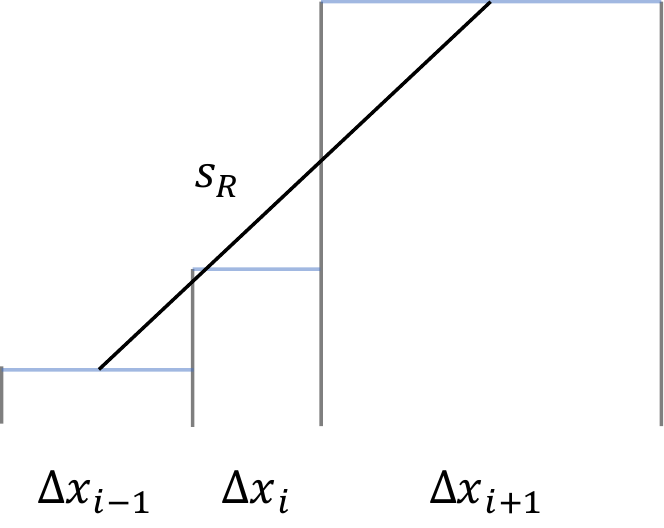}
    \caption{}
    \label{fig:reference-slope-irregular-mesh}
  \end{subfigure}
  \begin{subfigure}{0.45\textwidth}
    \centering
    \includegraphics[scale=0.8]{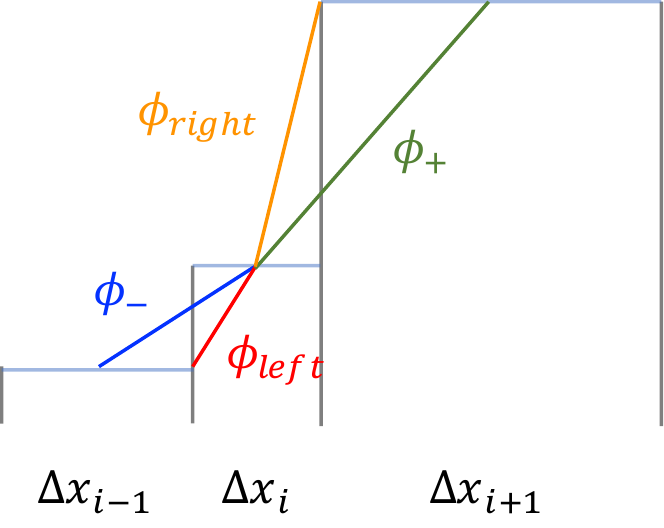}
    \caption{}
    \label{fig:s-minus-s-plus-irregular-mesh}
  \end{subfigure}

  \caption{(a) The reference slope on a irregular mesh;
    (b) Differences slopes on a irregular mesh.}
\end{figure}

Same as the uniform mesh, the reference slope $s_R$ is defined as the slope of the line that connects
($x_{i-1}, U_{i-1}$) and ($x_{i+1}, U_{i+1}$):
\begin{equation}
  s_R \coloneqq \frac{U_{i+1} - U_{i-1}}{x_{i+1} - x_{i-1}} = \frac{2 \Delta_t}{(2 + a + b) \Delta x_i},
\end{equation}
as illustrated in figure \ref{fig:reference-slope-irregular-mesh}.
The four bounding lines are shown in figure \ref{fig:s-minus-s-plus-irregular-mesh}, and their non-dimensional slopes are:
\begin{equation}
\begin{split}
  &{\color{blue}\phi_{-}} \coloneqq \frac{s_{-}}{s_R} = \frac{2+a+b}{1+a} f \\
  &{\color{red}\phi_{left}} \coloneqq \frac{s_{left}}{s_R} = (2+a+b) f \\
  &{\color{OliveGreen}\phi_{+}} \coloneqq \frac{s_{+}}{s_R} = \frac{2+a+b}{1+b} (1-f) \\
  &{\color{orange}\phi_{right}} \coloneqq \frac{s_{right}}{s_R} = (2+a+b) (1-f)
\end{split}
\end{equation}
in which, $f$ is the same as defined for uniform-sized mesh, and $f \coloneqq \Delta_{-}/\Delta_t$.
Similar to the uniform mesh case, at several special $f$ values, two of these four slopes are equal.
These special $f$ values are given as follows:
\begin{equation}
\label{eqn:special_f_irregular_mesh}
\begin{aligned}
  &f_1 = \frac{1}{2+b}      &\text{when } {\color{red}\phi_{left}} = {\color{OliveGreen}\phi_{+}}     &{}\\
  &f_2 = \frac{1+a}{2+a+b}  &\text{when } {\color{blue}\phi_{-}} = {\color{OliveGreen}\phi_{+}}       &{}\\
  &f_3 = \frac{1+a}{2+a}    &\text{when } {\color{blue}\phi_{-}} = {\color{orange}\phi_{right}}       &{}
\end{aligned}
\end{equation}

\begin{figure}[h]
  \centering
  \begin{subfigure}[b]{0.45\textwidth}
    \centering
    \includegraphics[scale=0.6]{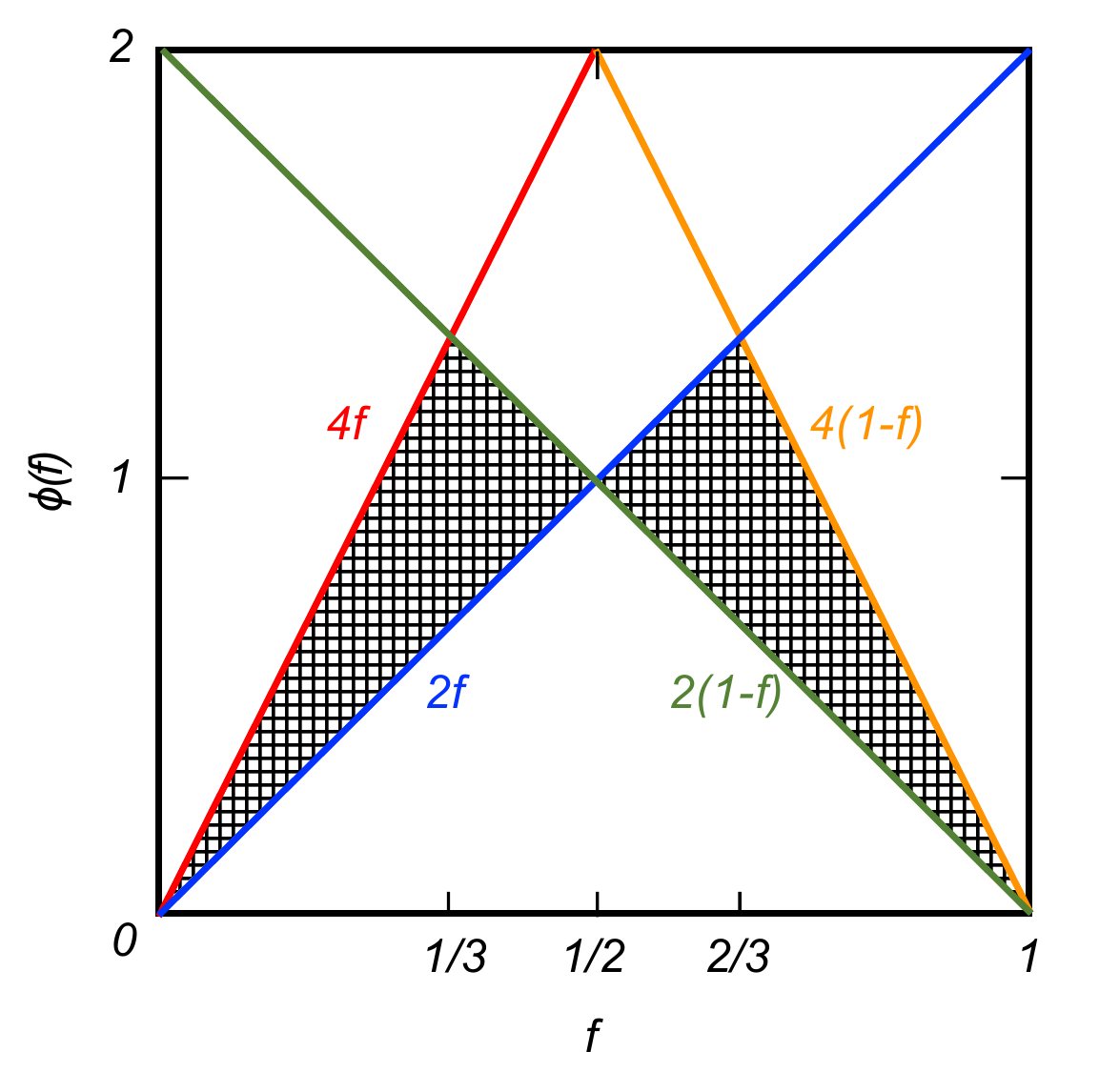}
    \caption{}
    \label{fig:uniform-mesh-tvd}
  \end{subfigure}
  \begin{subfigure}[b]{0.45\textwidth}
    \centering
    \includegraphics[scale=0.6]{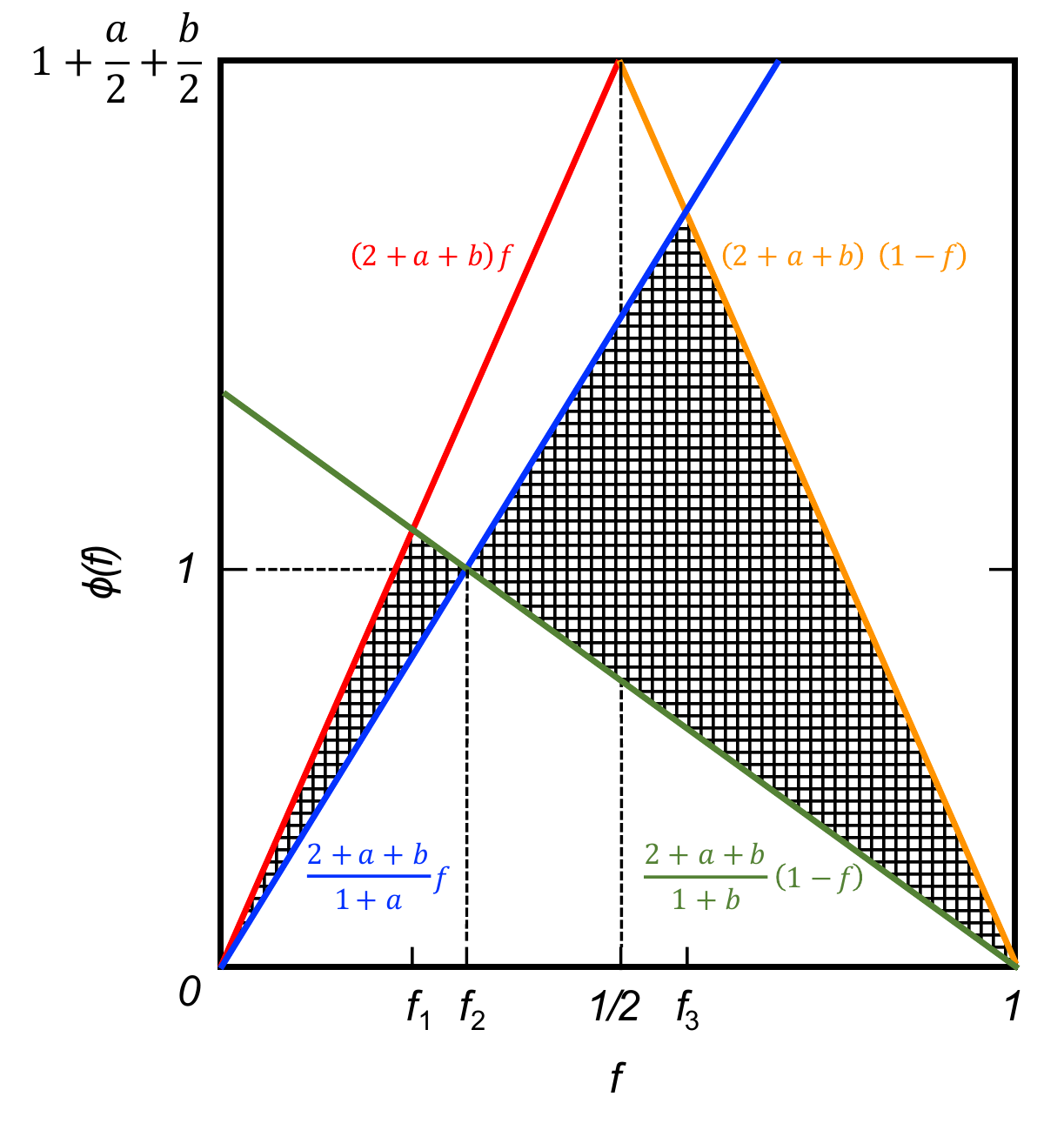}
    \caption{}
    \label{fig:irregular-mesh-tvd}
  \end{subfigure}

  \caption{Illustration of (a) the TVD region (shaded); and (b) the high-resolution TVD region (shaded) on the $\phi$-$f$ plot.
    It is clear that the high-resolution TVD region is a subset of the TVD region.}
  \label{fig:tvd-comparison}
\end{figure}

The direct outcome of \textbf{Conjecture 1} is that, on non-uniform meshes, the high-resolution TVD slope is restricted by
the two smallest slopes among the four colored lines.
On the $\phi$-$f$ plot, this is the shaded region illustrated in figure \ref{fig:irregular-mesh-tvd}.
In the same figure, the high-resolution TVD region on uniform mesh is also shown for comparison.
It is easy to see that figure \ref{fig:irregular-mesh-tvd} reduces to figure \ref{fig:uniform-mesh-tvd} on a uniform mesh, where $a = b = 1$.

At this stage, it is tempting to obtain high-resolution TVD scheme on non-uniform meshes by projecting
existing high-resolution TVD schemes on uniform meshes.
This can be done by finding the transform matrix that will homogeneously transform the two triangles of figure \ref{fig:uniform-mesh-tvd} into the two in \ref{fig:irregular-mesh-tvd}.
The transformed high-resolution TVD schemes can then be obtained by using this transform matrix.
This however does not seem to be the most efficient way, and there does not seem to have direct benefit to do so,
e.g., it is difficult to define symmetry on the irregular mesh, and the transformation of symmetric slope limiter does not seem to be very useful.

It is noted that figure \ref{fig:irregular-mesh-tvd} is equivalent to figure 4 of Berger et al. \cite{Berger_2005},
except that 1) Berger et al. used a different reference slope, so in our plot, the y-axis is scaled by a factor of $(2+a+b)/2$;
and 2) Berger et al. only marked the TVD region, not the high-resolution TVD region as we illustrate in figure \ref{fig:irregular-mesh-tvd}.
Berger et al. \cite{Berger_2005} also proposed two functions that lay in the TVD region,
which can be seen as generalized van Leer slope limiters,
as they recover the usual van Leer slope limiter on uniform meshes.
It can also be easily proved that, with proper scaling of $(2+a+b)/2$, they lay in the high-resolution TVD region proposed in this article.
To clarify the discussion, one of the slope limiter provided by Berger et al. \cite{Berger_2005}, i.e., their equation (38), is scaled and provided as:
\begin{equation}
  \phi(f) = \frac{2+a+b}{2} \left\{
  \begin{aligned}
    &f\left[1-\frac{a}{1+a}\left(\frac{f}{f_2}\right)^{1/a}\right]        & f \leq f_2 \\
    &(1-f)\left[1-\frac{b}{1+b}\left(\frac{1-f}{f_2}\right)^{1/b}\right]  & f > f_2
  \end{aligned}
  \right.
\end{equation}
As discovered by Berger et al. \cite{Berger_2005}, this slope limiter does a very good job on various non-uniform meshes.
The details on the performance of this slope limiter on non-uniform meshes are referred to Berger et al. \cite{Berger_2005}.
The same slope limiter has also been used by the author in a previous numerical simulation study of the Welander natural circulation problem, which showed very good results \cite{Zou_PNE}.

\section{Summary}
In this paper, we have illustrated how the concept of slope limiter can be interpreted graphically for uniform meshes.
We then conjecture that the same graphical rule can be generalized for non-uniform meshes.
The high-resolution total variance diminishing (TVD) region of slope limiter for non-uniform meshes can then be obtained.

Future work is needed to prove the conjecture, which seems possible if following Sweby's derivations for uniform meshes \cite{Sweby_1984}.

\section*{Acknowledgment}
Argonne National Laboratory’s work was supported by the U.S. Department of Energy, Office of Nuclear Energy, under contract DE-AC02-06CH11357.


\end{document}